\documentclass[a4paper,11pt]{article}

\usepackage{cite}
\usepackage[pdftex]{graphicx}
\usepackage[cmex10]{amsmath}
\usepackage{amsfonts} 

\begin{document}

\title{A New Mathematical Model for Evolutionary Games on Finite Networks of Players}
\author{Dario~Madeo and
        Chiara~Mocenni
\thanks{D. Madeo and C. Mocenni are  with the Department
of Information Engineering and Mathematical Science, University of Siena, 
Via Roma 56, 53100 - Siena (Italy). E-mail: \{madeo,mocenni\}@dii.unisi.it}.
}
\date{July 3, 2013}

\maketitle

\begin{abstract}
A new mathematical model for evolutionary games on graphs is proposed to 
extend the classical replicator equation to finite populations of players 
organized on a network with generic topology. Classical results from game theory, 
evolutionary game theory and graph theory are used. More specifically, 
each player is placed in a vertex of the graph and he is seen as an 
infinite population of replicators which replicate within the vertex. 
At each time instant, a game is played by two replicators belonging to 
different connected vertices, and the outcome of the game influences 
their ability of producing offspring. Then, the behavior of a 
vertex player is determined by the distribution of strategies used by
the internal replicators.  Under suitable 
hypotheses, the proposed model is equivalent to the classical replicator 
equation. Extended simulations are performed to show the dynamical 
behavior of the solutions and the potentialities of the developed model.
\end{abstract}


\section{Introduction}

In mathematical literature, there are important examples of models 
designed to describe the dynamical interaction between a set of players 
in a game-like context \cite{HS2003,N2006,NM1944,HS1988}. Players are 
undistinguishable members of a large population, characterized by a 
phenotype which determines the fixed strategy they choose among the $M$ 
available, when playing with any other individual randomly selected in 
the population.\\

The payoff earned in the game by each player depends on the elements of 
a game-specific payoff matrix, while the system dynamics is described by 
an ordinary differential equation defined in the N-simplex, namely the 
replicator equation \cite{HS1998,MS1982}.\\ The solutions of the replicator 
equation may evolve towards evolutionary stable strategies, which are 
Nash equilibria of the game and asymptotically stable stationary states 
\cite{N1950}. Evolutionary stable strategies are robust against invasion 
by competing strategies \cite{NSTF2004, MS1982}. Other kind of Nash 
equilibria can also exist, such as equilibria corresponding to Lyapunov 
stable stationary states of the replicator equation \cite{W1995}. 
\\

Replicator equation has been used to describe several phenomena, such as biological 
evolution driven by replication and selection \cite{NS2004}, and reaction-diffusion 
dynamics \cite{KE2000,SPB2011}. Replicator dynamics including mutation have been studied 
in \cite{N2006} and may lead to more complex behaviors, characterized by
 Hopf bifurcation and limit cycles \cite{PL2011,PCL2012}.
Moreover, the replicator equation has been used for solving decision and consensus control 
problems \cite{EY2010,BGP2008} and machine learning for  optimization \cite{P1999}.
Significant application have been also developed in the field of social
science \cite{BJR2003}.
Although widespread, the replicator equation is grounded on several 
strong assumptions on the system under investigation. 

\begin{enumerate}
\item The population is very large. 
Indeed, dominating strategies emerge due to higher rates of replication 
than the others, causing the frequency of less efficient strategies to 
become irrelevant in the total population.
\item Any member of the population can play with any other member with 
the same probability. Participants of each game are chosen randomly and 
no social structures are present in the population.
\item The payoff earned by each player is defined by the payoff matrix, 
which is unique for all the population. In some cases, two or more
subpopulations with different payoff matrices are considered \cite{SF2007}. 
\item Players are constrained to behave according to a single fixed 
strategy at each round of the game they are playing. 
For instance, if the phenotype of an individual is to be aggressive or generous, 
he will show the same level of aggressiveness or generosity in any 
situation at any time, without having the capacity of regulating the 
level of his natural impulses. 
\end{enumerate}

The above assumptions are the basis of the well known equivalence between
the replicator equation and some ecological models, such as the 
predator-prey model introduced by Lotka and Volterra \cite{HS1998}.\\

Several modifications have been introduced in the replicator equation 
to overcome the limitations caused by the above assumptions. 
For example, in \cite{BJR2003} is presented a
generalization of the replicator equation to $N$ players.
Concerning networked populations, 
many efforts have been done to include the topology of connections between 
players  in \cite{GR2012,SRP2006,N2006,TK2002,TFSN2004} to deal with scenarios 
in which the connection between agents plays a fundamental role and may 
yield to the interpretation of inspected phenomena. For example, in 
\cite{GR2012} is proposed an algorithm that uses the connections among 
players and specific updating rules to induce cooperative behavior in 
an evolutionary prisoner dilemma game. On the other hand, a seminal paper by Ohtsuki et al. \cite{ON2006} 
presented the replicator equation on infinite graphs, under the assumption 
that every vertex is connected to the same number of neighbours. In 
particular, he showed that the replicator equation on a graph with 
fixed degree is equivalent to a classical replicator equation when 
the degree goes to infinity. \\

In this paper, we derive a suitable mathematical model to describe static and
dynamical behaviour of an $N$-players game interaction, where an agent 
is meant to engage his challenges only with a restricted set of players
connected to him. To this aim, we 
make use of the standard non-cooperative game theory results and 
graph theory.
The results presented in the paper generalize many assumptions under the 
classical replicator equation and recent results on the replicator 
equation on graphs. Specifically, the equation derived in this paper 
meets the following:
\begin{enumerate} 
\item A finite (even small) or infinite population is considered.
\item The elements of the population are the vertices of a graph. 
Each element can be engaged in a game with another element only if they 
are connected in the graph. No constraints on the topology of the graph 
are assumed. Moreover, the connections can be weighted to remark 
different perceived importance of each interaction.
\item The payoff matrix can be player-specific, including the situations
where the perception of the game is different for each individual.
\item Each player can behave according to a combination of strategies. 
He is a sort of "mixed player", thus incorporating composite and 
multiple personality traits. His behavior can be driven contemporarily 
by heterogeneous impulses with different strengths, such as, for 
example, being cooperative and non cooperative, generous and selfish, 
at the same time. The proposed approach makes players more realistic than 
in the classical framework and naturally extends the evolutionary game 
theory to a social context with human players.
\end{enumerate}

The new framework presented in the paper generalizes the classical 
replicator equation, that  can be obtained as a special case by assuming 
that any individual of the population possesses the same payoff matrix 
and starts playing from an identical initial condition.\\

The paper is structured as follows. Section \ref{sec:ncgames} presents 
some preliminaries on noncooperative games on graphs. Then, the extended version 
of the replicator equation on graphs with generic topology is introduced 
in Section \ref{sec:reg}. Some properties of the new replicator equation 
are presented in Section \ref{sec:mathreg}, including 
the equivalence to the classical replicator equation when homogeneous 
initial conditions are used. Extended simulations are reported in 
Section \ref{sec:sim}, while some conclusions and future work are 
discussed in Section \ref{sec:conclusions}.

 
\section{Non-cooperative games on graphs}\label{sec:ncgames}

In real world situations, interactions between a finite number of rational
players can be influenced by topological constraints; in most cases, each
player is only able to meet with a reduced number of opponents which
are close with respect to a suitable topology, such as the distance
between them. In this sense, we talk about networks of players, where
interconnections depend on the context. An interesting case of interaction
is represented by non-cooperative games, extensively studied in 
\cite{NM1944,N1950}.\\

In this section, we extend the classical game theory by introducing a 
network, represented by a graph, which describes the connection among the
involved players.

\subsection{Preliminaries about games on graphs} 
\label{sec:nmgames}

Typically, networks are described
by means of graphs, and in a game context, each player is represented by a vertex.
An edge between two players indicates that they interact. However, a player
can consider that some interactions are more important than others.
Moreover, two connected players can have different perception 
of the importance of their interaction. These aspects can
be accounted by assuming that the graph is weighted and directed;
an edge starting from a player and ending to another, is labeled with
a positive weight to indicate the importance that the first player
attributes to the game.\\

Formally, let $\mathcal{G}$ be a directed weighted graph of order $N < +\infty$, 
and let $\mathcal{V}$ be the set of vertices (players) ($\mathcal{V} = \{1, \ldots,N\}$).
The graph $\mathcal{G}$ is fully described by its adjacency 
matrix $\boldsymbol{A} \in \mathbb{R}_{+}^{N \times N}$; in particular, when player $v$
is meant to play with player $w$, then there is an edge which starts from
$v$ and ends to $w$. In this case, $(v,w)$-entry of $\boldsymbol{A}$, $a_{v,w}$, 
is the positive weight attributed by $v$ to
the game against $w$. When $a_{v,w} > 0$ and $a_{w,v} = 0$, there is an 
interaction between $v$ and $w$, but only $v$ will get a payoff after the
challenge. Finally, if both $a_{v, w}$ and $a_{w,v}$ are equal to $0$, then there is no
interaction between these players. In general, $a_{v, w} \neq a_{w, v}$.
We assume that $\mathcal{G}$ has no self-edge, which means that no player 
has interaction with himself (i.e. $a_{v,v} = 0 ~~ \forall v \in \mathcal{V}$). 
We indicate with $\mathcal{N}_{v}  = \{w \in \mathcal{V} : a_{v,w} > 0 \vee a_{w,v} > 0\}$ the
neighborhood of $v$ (i.e. the set of vertices that interact with $v$), and
with $\mathcal{N}^{+}_{v}  = \{w \in \mathcal{V} : a_{v,w} > 0\}$ the
out-neighborhood of $v$ (i.e. the set of vertices $v$ is connected to with
and exiting edge). The cardinalities of these set are indicated with
$\delta_{v}$ and $\delta^{+}_{v}$, and they represent the degree and out-degree
of $v$ respectively. Note that, in general, $\mathcal{N}^{+}_{v} \subseteq 
\mathcal{N}_{v}$, and also $\delta^{+}_{v} \leq \delta_{v}$.\\

A player $v$ will play exactly $\delta_{v}$ two-players games with all 
its neighbours. In each one-to-one competition, the set of available 
strategies for both players is $\mathcal{S}=\{1, \ldots, M\}$, while the 
outcome that player $v$ can obtain is defined by 
a payoff matrix $\boldsymbol{B}_{v} \in \mathbb{R}^{M \times M}$; 
when player $v$ uses 
strategy $s \in \mathcal{S}$ in a two-players game against a player which
uses strategy $r \in \mathcal{S}$, then he earns a payoff equal to the 
$(s,r)$-entry of $\boldsymbol{B}_{v}$, $b_{v,s,r} = \boldsymbol{e}_{s}^{T}\boldsymbol{B}_{v}\boldsymbol{e}_{r}$, where
$\boldsymbol{e}_{s}$ and $\boldsymbol{e}_{r}$
are the $s$-th and $r$-th versors of $\mathbb{R}^{M}$, respectively. \\

Each player decides to use the same strategy $s \in \mathcal{S}$
in all the games he is involved in. He will play against all
vertices in $\mathcal{N}_{v}$, but he will earn a payoff only when he plays
with a player $w \in \mathcal{N}^{+}_{v}$, since when $a_{v, w}=0$ and $a_{w,v}>0$, 
there is an interaction which is meaningful only for player $w$.\\

\subsection{Effective payoff for games on graphs}\label{sec:payoff_graphs}

In an interconnected context, the effective payoff earned 
(or the \textit{fitness} of a strategy) must be defined as an \textit
{environmental measure} depending on all the 
interactions between near players. This measure must quantify 
how well a strategy behaves. Since each connection between two players
has a positive weight, we pose that the effective payoff for a generic
player $v$ is the weighted average of all obtained payoffs.
Let's denote with $s_{w} \in \mathcal{S}$ the strategy of the generic 
player $w$. Then, the effective payoff of player $v$, $\pi_{v}(s_{1}, \ldots, s_{N})$ is the following:
\begin{eqnarray}
\displaystyle \pi_{v}(s_{1}, \ldots, s_{N}) & = & \displaystyle  \frac{1}{d_{v}} \sum_{w \in \mathcal{N}^{+}_{v}} \boldsymbol{e}_{s_{v}}^{T} \boldsymbol{B}_{v} \boldsymbol{e}_{s_{w}} = \nonumber \\
& = & \displaystyle  \boldsymbol{e}_{s_{v}}^{T} \boldsymbol{B}_{v} \left( \frac{1}{d_{v}} \sum_{w = 1}^{N} a_{v, w}  \boldsymbol{e}_{s_{w}} \right), \label{eqn:payoff_pure_wa}
\end{eqnarray}
where $d_{v} = \sum_{w = 1}^{N}a_{v,w}$ is the normalization factor. This
model of payoff based on weighted average will be denoted with WA.
However, there are situations
in which payoffs are cumulative and the weighted sum is used without 
the normalization factor $d_{v}$. In this case we have that:
\begin{eqnarray}
\displaystyle \pi_{v}(s_{1}, \ldots, s_{N}) & = & \displaystyle  \sum_{w \in \mathcal{N}^{+}_{v}} \boldsymbol{e}_{s_{v}}^{T} \boldsymbol{B}_{v} \boldsymbol{e}_{s_{w}} = \nonumber \\
& = & \displaystyle  \boldsymbol{e}_{s_{v}}^{T} (d_{v}\boldsymbol{B}_{v}) \left( \frac{1}{d_{v}} \sum_{w = 1}^{N} a_{v, w}  \boldsymbol{e}_{s_{w}} \right). \label{eqn:payoff_pure_ws}
\end{eqnarray}
The payoff model based on weighted sum (WS) can be considered as 
WA, where each payoff matrix is substituted by $d_{v}\boldsymbol{B}_{v}$. For this reason,
we will mainly work on WA model, unless differently specified. \\

The term $\frac{1}{d_{v}} \sum_{w = 1}^{N} a_{v, w}  \boldsymbol{e}_{s_{w}} $
that appears in both WA and WS models, is a vector where all components
are non-negative numbers which sum up to $1$. In a certain way, 
player $v$ fights against one \textit{virtual player} which summarize all 
the strategies
used by its opponents in the set $\mathcal{N}^{+}_{v}$; in general, the
strategy used by the virtual player is a mixed strategy which
represents what player $v$ effectively \textit{sees} around him. This aspect
will be deeply investigated later in this paper, because it plays a 
fundamental role to reach our aim. \\

\subsection{$(N-M)$-games and games on graphs}
\label{sec:nmgames}

Notice that, for each $v$, $\pi_{v}$ can be interpreted as a $N$-dimensional 
tensor, where the $(s_{1}, \ldots, s_{N})$-entry is $\pi_{v}(s_{1}, \ldots, s_{N})$.
In this way, the game interaction between
interconnected players on a finite graph is equivalent to a $N$-players game,
where the set of pure strategies is $\mathcal{S}$, and the payoff
of player $v$ is represented by the tensor $\pi_{v}$. The structure of the
graph is embedded in this definition, since the payoff tensor 
depends on the adjacency matrix $\boldsymbol{A}$. Moreover, there are no assumptions made on
the structure of the graph itself.\\

For example, consider the following matrices:
\begin{equation}\label{eqn:ex_AandB}
\displaystyle \boldsymbol{A} = \left[
\begin{array}{ccc}
0 & 1 & \mu \\
\mu & 0 & 2\mu \\
0 & \mu & 0
\end{array}
\right], ~~\boldsymbol{B} = \left[{}
\begin{array}{cc}
a & b \\
c & d 
\end{array}
\right],
\end{equation}
where $\mu \geq 0$, and assume that $\boldsymbol{B}_{v} = \boldsymbol{B}$ 
for all $v$.
In this case, $\mathcal{V} = \{1, 2, 3\}$ and $\mathcal{S} = \{1, 2\}$.
Table \ref{table:payoffs} shows the payoff tensors $\pi_{v}$ of each player, 
which depend on the model parameters $\mu, a, b, c$ and $d$. Both models WA and WS are considered. \\

\begin{table}
\begin{center}
\begin{tabular}{| c | c | c || c | c || c | c || c | c |}
\hline
\multicolumn{3}{|c||}{} & \multicolumn{2}{|c||}{$\pi_{1}(s_{1}, s_{2}, s_{3})$} & \multicolumn{2}{|c||}{$\pi_{2}(s_{1}, s_{2}, s_{3})$} & \multicolumn{2}{|c|}{$\pi_{3}(s_{1}, s_{2}, s_{3})$} \\
\hline
$s_{1}$ & $s_{2}$ & $s_{3}$ & WA & WS & WA & WS & WA & WS \\ \hline\hline
$1$ & $1$ & $1$ & $a$ & $a(1+\mu)$ & $a$ & $3a\mu$ & $a$ & $a\mu$ \\ \hline
$1$ & $1$ & $2$ & $\frac{a+ b\mu}{1 + \mu}$ & $a+ b\mu$ & $\frac{a+ 2b}{3}$ & $(a+2b)\mu$ & $c$ & $c\mu$ \\ \hline
$1$ & $2$ & $1$ & $\frac{a\mu+ b}{1 + \mu}$ & $a\mu + b$ & $c$ & $3c\mu$ & $b$ & $b\mu$ \\ \hline
$1$ & $2$ & $2$ & $b$ & $b(1 + \mu)$ & $\frac{c+ 2d}{3}$ & $(c+2d)\mu$ & $d$ & $d\mu$ \\ \hline
$2$ & $1$ & $1$ & $c$ & $c(1+\mu)$ & $\frac{2a+ b}{3}$ & $(2a+b)\mu$ & $a$ & $a\mu$ \\ \hline
$2$ & $1$ & $2$ & $\frac{c+ d\mu}{1 + \mu}$ & $c+d\mu$ & $b$ & $3b\mu$ & $c$ & $c\mu$ \\ \hline
$2$ & $2$ & $1$ & $\frac{c\mu+ d}{1 + \mu}$ & $c\mu+d$ & $\frac{2c+ d}{3}$ & $(2c+d)\mu$ & $b$ & $b\mu$ \\ \hline
$2$ & $2$ & $2$ & $d$ & $d(1+\mu)$ & $d$ & $3d\mu$ & $d$ & $d\mu$ \\ \hline
\end{tabular}
\end{center}
\caption{Payoff tensor of the game on graphs defined by matrices in equation \eqref{eqn:ex_AandB}
for  WA and WS payoff models. For each combination of strategies ($s_{1}, s_{2}, s_{3}$) the payoffs $\pi_{1}, \pi_{2}$ and $\pi_{3}$
of players $1$, $2$, $3$ are reported.}
\label{table:payoffs}
\end{table}

It is evident that the presence of weights, the asymmetry of the matrix, and
the use of a particular payoff model may lead to very different calculation
of the payoff tensor, and hence, the structure of the game itself changes.
Indeed, the effective payoff obtained by 
a player when he is engaged in a game is essentially evaluated by means 
of tensors, depending on the adjacency matrix of the graph.
These payoffs  define the \textit{virtual player} mentioned at the end 
of Section \ref{sec:payoff_graphs}, which embodies all the strategies
used by the player's opponents. As a consequence, each player in the 
game is a sort of "mixed player", thus incorporating composite and 
multiple personality traits, and behaving according to heterogeneous 
impulses with different strengths. \\

As a natural consequence, a $N$-players $M$-strategies game (from now on, 
$(N,M)$-game) can be extended over the set of mixed strategies $\Delta_{M}$: 
\begin{equation}\nonumber
\displaystyle\Delta_{M} = \{\boldsymbol{z}=[z_1 \ldots z_{M}]^{T} \in \mathbb{R}^{M} : \sum_{i=1}^{M} z_{i} = 1 \wedge z_{i} \geq 0 ~~ \forall i \in \mathcal{S} \}.
\end{equation}

We indicate with $\boldsymbol{x}_{v} = [x_{v,1} \ldots x_{v,M}]^{T} \in \Delta_{M}$ the mixed
strategy of player $v$. Recall that $x_{v,s}$ is the probability that
player $v$ uses strategy $s$, while he takes part in the games. The formula of the 
expected effective payoff that player $v$ obtains, is similar to equation \eqref{eqn:payoff_pure_wa}:
\begin{eqnarray}
\displaystyle \pi_{v}(\boldsymbol{x}_{1}, \ldots, \boldsymbol{x}_{N}) & = & \displaystyle  \frac{1}{d_{v}} \sum_{w \in \mathcal{N}^{+}_{v}} \boldsymbol{x}_{v}^{T} \boldsymbol{B}_{v} \boldsymbol{x}_{w} = \nonumber \\
& = & \displaystyle \boldsymbol{x}_{v}^{T} \boldsymbol{B}_{v} \left( \frac{1}{d_{v}} \sum_{w = 1}^{N} a_{v, w}  \boldsymbol{x}_{w} \right), \label{eqn:payoff_mixed}
\end{eqnarray}
where $\boldsymbol{x}_{v}^{T} \boldsymbol{B}_{v} \boldsymbol{x}_{w}$ represents the expected outcome for player $v$ of the one-to-one game played
by $v$ itself against $w$. From now on, we pose that 
$\pi_{v}(\boldsymbol{x}_{v}, \boldsymbol{x}_{-v})=\pi_{v}(\boldsymbol{x}_{1}, \ldots, \boldsymbol{x}_{N})$, where
$\boldsymbol{x}_{-v}$ indicates the group of all the vectors $\boldsymbol{x}_{w}$, with $w \neq v$.\\

When a pure strategy $s$ is used by player $v$, vector $\boldsymbol{x}_{v}$
takes the form of the standard versor $\boldsymbol{e}_{s}$ of $\mathbb{R}^{M}$.
If each $\boldsymbol{x}_{w}$ is a versor, say $\boldsymbol{x}_{w} = \boldsymbol{e}_{s_{w}}$, 
then equations \eqref{eqn:payoff_pure_wa}
and \eqref{eqn:payoff_mixed} coincide. Furthermore, we can define 
the expected payoff $p_{v,s}$ of
player $v$ when he is preprogrammed to use the strategy $s$ 
(hence $\boldsymbol{x}_{v} = \boldsymbol{e}_{s}$) in all games played against its neighbours:
\begin{eqnarray}
\displaystyle p_{v,s} & = & \displaystyle \pi_{v}(\boldsymbol{e}_{s}, \boldsymbol{x}_{-v}) = \nonumber\\
& = & \displaystyle  \boldsymbol{e}_{s}^{T}\boldsymbol{B}_{v} \left(\frac{1}{d_{v}}\sum_{w=1}^{N}a_{vw}\boldsymbol{x}_{w}\right). \label{eqn:pvs}
\end{eqnarray}
Equation \eqref{eqn:pvs} easily leads to a more convenient definition of the expected payoff 
obtained by the player $v$. That is:
\begin{eqnarray}
\displaystyle \phi_{v} & = & \displaystyle\pi_{v}(\boldsymbol{x}_{v}, \boldsymbol{x}_{-v}) = \nonumber \\
& = & \displaystyle \boldsymbol{x}_{v}^{T}\boldsymbol{B}_{v} \left(\frac{1}{d_{v}}\sum_{w=1}^{N}a_{vw}\boldsymbol{x}_{w}\right). \label{eqn:phiv}
\end{eqnarray}

The present work uses the same theoretical issues developed in the classical 
non-cooperative $(N-M)$-games theory (\cite{NM1944, N1950, BCV1997}). 
The differences here introduced, consist with the possibility of embedding
any topological structure in the game; indeed, the payoff tensor used to describe
the game depends on the adjacency matrix of the graph. 

 
\section{The replicator equation on graphs}\label{sec:reg}

Thanks to game theory, we are able to predict the strategies of opponents, 
assuming that all of them behave in a rational way during their
decision-making tasks. In fact, rational players choose pure strategies 
which may lead to a pure Nash equilibrium, whenever it exists. 
Recall that Nash theorem \cite{N1950} asserts that a game has always at 
least one Nash equilibrium within the set of mixed strategies.
If the game is repeated over time,
we can imagine that a mixed strategy describes, for each time, the 
probability that a player uses a certain pure strategy independently from 
the choices made until that moment. 
Payoff is computed as the average outcome that a player obtains when the
game is reiterated for an infinite number of times, and Nash equilibria 
are evaluated accordingly.
Although players' behavior is randomized over time, there is a precise 
rational scheme that they follow. Hence, pure and mixed 
strategies games are quite similar, since in both cases players decide
their behaviors through a rational decision-making task at the beginning of the game. \\

Often, in real world situations, players do not have a full knowledge 
about the game, and the decision-making process suggested by game
theory is not applicable. However, players can learn from the context; 
time after time, they are able to compare their payoffs with the 
outcomes of their opponents, and strategies are changed accordingly. 
For example, after some game iterations, one can understand that
its opponent is preprogrammed to play always the same strategy $r$, and 
then he decides to adopt the strategy $s$ which is the best reply to $r$
(i.e. he obtains the maximum payoff knowing the strategy of his opponent). Also,
a player can simply imitate the strongest opponent to reach a 
greater payoff.\\

Evolutionary biology gives an interesting interpretation for mixed
strategies. Let's consider a large population in which each individual is
programmed to play a particular pure strategy. Population is divided 
according to a mixed strategy, which indicates the frequency of 
pure strategies in the population.  
Pairs of individuals are randomly drawn from the population to play games. 
The average payoff obtained by all players playing a certain strategy 
is a measurement of the \textit{fitness} of that strategy. \\

Nature promotes fittest strategies; time after time, when a strategy 
has a payoff greater than average, its frequency must increase, 
and consequently, frequencies relative to poorly fit strategies decrease. 
This \textit{natural selection} process may yield to interesting
dynamical phenomena, since the fitness of each strategy changes over
time according to the frequencies of the subpopulations.
Natural selection is realized through reproduction. 
Generally, this happens asexually and the offspring produced is identical
to parent. For these reasons, a player in this context is also known
as \textit{replicator}. He is preprogrammed to use a certain strategy and 
there is no rational decision-making process. A strategy is something
included in the genes of a replicator, and hence, exhibited from birth. 
A replicator uses a 
``good'' strategy when has a fitness greater than the average. 
Nature favors players with a good strategy, by allowing them 
to survive and to produce offspring.\\

The dynamics of replicator
populations can be used to describe several non-biological situations: 
scientific ideas, life-styles, political orientations diffuse by means of 
imitation and education process that easily replace the concept of asexual
reproduction. In fact, this becomes clearer if we interpret fitness 
as a measurement quantifying how well a strategy behaves in a certain 
context, whatever it is. 

\subsection{Towards evolutionary $(N-M)$-games}\label{sec:towardsNMgames}

In this section we develop a key idea which leads to the definition of
a replicator equation based on a generic graph, where the 
number of players and the network structure are arbitrary.\\

\textbf{Main idea} - We imagine that each vertex of the network contains 
an infinite  population of individuals. 
We will refer to such elements as 
\textbf{atomic players}, and to vertices of the network as 
\textbf{vertex players}. The first are replicators, while the latter 
are the players introduced in paragraph \ref{sec:nmgames} for $(N-M)$-games.
Basically, each atomic player behaves like the corresponding vertex player:
an atomic player of $v$ takes part to $2$-players games, described by payoff matrix $B^{v}$, 
against exactly one atomic player randomly drawn from each connected
vertex, and his effective payoff is the average of
the payoffs obtained in all the one-to-one competitions.
By the way, atomic players are different from vertex players.
Indeed, all atomic players inside a vertex are indistinguishable,
except for the fact that each of them is preprogrammed to use a certain 
pure strategy in $\mathcal{S}$ during all the games he is involved in. 
On the contrary, vertex players can also adopt mixed strategies. 
Atomic players reproduce themselves by replication, after their participation 
to games, inside 
their population. Furthermore, their capacity to produce offspring is related 
to the effective payoff obtained.\\

Assume that one atomic player is randomly draw from each populations,
and let $s_{w}$ be the strategy used by the one extracted from population $w$. 
Then, the effective payoff earned by the atomic player of population $v$ is
defined by equation \eqref{eqn:payoff_pure_wa}.
Now, $x_{v,s}$ can be interpreted as the share of atomic players 
preprogrammed to use the pure strategy $s$ inside the vertex $v$. This implies
that $p_{v,s}$  is the expected effective payoff 
obtained by an atomic player of $v$ when he uses strategy $s$, while
$\phi_{v}$ represents the expected effective payoff for a generic atomic player
randomly drawn from population $v$ (see equations \eqref{eqn:pvs} 
and \eqref{eqn:phiv}).
A mixed strategy of a vertex player can be interpreted as the way in which 
its internal population of atomic players is distributed,
according to the pure strategy they are preprogrammed to play.

\subsection{Mathematical formulation of the replicator equation on graphs}\label{sec:mathreg}

Suppose now that the games are iterated in time. We will refer to game
session as the whole set of $2$-players games performed on the graph.
The probability for a replicator to survive and to reproduce himself
between two games sessions, depends on the comparison between the effective 
payoff obtained and the average effective payoff of all other players.
Let's assume that games' sessions take place at 
discrete and equidistant times (say, a session after each $\tau$ seconds). 
Let $x_{v,s}(t)$ indicate the share of population inside vertex $v$, which is
preprogrammed to use the pure strategy $s$, at time $t$.\\

What happens to $x_{v,s}(t + \tau)$? First of all, suppose for a while that the
population size in the generic node $v$ at time $t$ is $n_{v}(t)$. 
According to \cite{MS1982,HS2003,N2006}, we can consider that 
$p_{v,s}(t)$ represents a reproductive rate and therefore, $p_{v,s}(t)\tau $ 
is the number of offspring produced by one atomic player in $v$ that uses 
strategy $s$ between $t$ and $t+\tau$. Hence, the population size after 
a time $\tau$ is equal to the previous size plus the produced offspring,
since each atomic player reproduces himself within his population. That is:
\begin{equation}\nonumber
\displaystyle n_{v}(t + \tau) = n_{v}(t) + \sum_{r=1}^{M} n_{v}(t)x_{v,r}(t)p_{v,r}(t)\tau,
\end{equation}
where $n_{v}(t)x_{v,r}(t)$ is the size of the subpopulation which uses 
strategy $r$ at time $t$, and $n_{v}(t)x_{v,r}(t)p_{v,r}(t)\tau$ is the 
number of offspring produced by this subpopulation. By definition, 
$x_{v,s}(t + \tau)$ is the ratio between the size of subpopulation $s$ 
and the total population. Therefore:
\begin{eqnarray}
\displaystyle x_{v,s}(t + \tau) & = & \displaystyle\frac{n_{v}(t)x_{v,s}(t) + n_{v}(t)x_{v,s}(t)p_{v,s}(t)\tau}{n_{v}(t+\tau)} =\nonumber \\
& = & \displaystyle\frac{n_{v}(t)x_{v,s}(t)(1  + p_{v,s}(t)\tau)}{n_{v}(t)\left(1 + \displaystyle\sum_{r=1}^{M} x_{v,r}(t)p_{v,r}(t)\tau\right)} =\nonumber\\
& = & \displaystyle\frac{x_{v,s}(t)(1  + p_{v,s}(t)\tau)}{1 + \phi_{v}(t)\tau}.\label{eqn:x_map}
\end{eqnarray}
Notice that equation \eqref{eqn:x_map} does not depend on $n_{v}(t)$, and
hence this relationship is valid for any starting size of the population.\\

Our aim is to develop a mathematical model that describes the evolutionary
process on a graph when the time between replication events goes 
to $0$, thus making atomic players able to reproduce themselves continuously
in time. Let's consider the difference ratio of $x_{v,s}(t)$: 
\begin{equation}
\displaystyle \frac{x_{v,s}(t + \tau) - x_{v,s}(t)}{\tau} = \displaystyle\frac{x_{v,s}(t)(p_{v,s}(t) - \phi_{v}(t))}{1 + \phi_{v}(t)\tau}.\nonumber
\end{equation}
Letting $\tau \rightarrow 0$, we obtain that:
\begin{equation}\label{eqn:ode}
\displaystyle \dot{x}_{v,s}(t)  = x_{v,s}(t)(p_{v,s}(t) - \phi_{v}(t)),
\end{equation}
where the ``dot'' indicates the derivative with respect to time $t$. Finally,
we can write the following Cauchy problem:
\begin{equation}\label{eqn:cauchy}
\displaystyle \left\{
\begin{array}{l}
\displaystyle \dot{x}_{v,s}(t) = x_{v,s}(t)(p_{v,s}(t) - \phi_{v}(t)) \\
\\
\displaystyle x_{v,s}(0) = c_{v,s}
\end{array}
\right. ~~ \forall v \in \mathcal{V}, ~~ \forall s \in \mathcal{S},
\end{equation}
where, for consistency, it is assumed that the distribution of strategies at the initial 
time $t=0$ is known for each vertex (i.e. $\boldsymbol{x}_{v}(0) = [c_{v,1} \ldots c_{v,M}]^{T} \in \Delta_{M}$).\\

Systems \eqref{eqn:cauchy} represents the replicator 
equation on a graph. Note that no assumptions on the structure of
the graph is needed to derive the equation \eqref{eqn:ode}. Indeed,
the adjacency matrix of the network is fully embedded in the payoff
tensors. \\

It is straightforward to note that the equation \eqref{eqn:ode} has 
a structure similar to the classical replicator equation; for example,
dominant strategies are the fittest, and hence when the relative fitness
$p_{v,s}$ is better than the average $\phi_{v}$, the corresponding frequencies
will grow over time. In the next section, the very strong correlation 
between the two equations will be rigorously shown. Furthermore, the relationship
between Nash equilibria of the underlying $(N-M)$-game and the rest points
of the dynamical equation \eqref{eqn:ode} will also be discussed in section
\ref{sec:reproperties}.


\section{Properties of the replicator equation on graphs}\label{sec:reproperties}

\subsection{Invariance of $\Delta_{M}$}

Let $\boldsymbol{x}_{v}(t)$ be the unique solution of problem \eqref{eqn:cauchy},
obtained by posing $\boldsymbol{x}_{v}(0) \in \Delta_{M}$. In addition, 
suppose that there exists a time instant $t_{2}$ where $x_{v,s}(t_{2}) < 0$.
Since all the components of the solution are continuous and non-negative at $t=0$, then
there must be a time $t_{1} < t_{2}$ such that $x_{v,s}(t_{1}) = 0$.
Following equation \eqref{eqn:ode}, we can state that $\dot{x}_{v,s}(t_1) = 0$, 
and hence, this component will be $0$ for all times after $t_1$. 
For the unicity of the solution, this implies that no time
$t_{2}$ for which $x_{v,s}(t_{2}) < 0$ exists. Thus, for each $v \in \mathcal{V}$
we have that:
\begin{equation}\label{eqn:nonnegativity}
\displaystyle \boldsymbol{x}_{v}(0) \in \Delta_{M} \Rightarrow x_{v,s}(t) \geq 0,
\end{equation}
for all strategies $s \in \mathcal{S}$ and for all times $t > 0$.
Notice that the total variation of the strategies distribution in a 
vertex is null at time $t$ when $\sum_{s=1}^{M}x_{v,s}(t) = 1$.
In fact:
\begin{eqnarray}
\displaystyle \sum_{s=1}^{M}\dot{x}_{v,s}(t) & = & \displaystyle\sum_{s=1}^{M}x_{v,s}(t)(p_{v,s}(t) - \phi_{v}(t)) = \nonumber\\
& = & \displaystyle\sum_{s=1}^{M}x_{v,s}(t)p_{v,s}(t) - \phi_{v}(t)\sum_{s=1}^{M}x_{v,s}(t) = \nonumber\\
& = & \displaystyle \phi_{v}(t) - \phi_{v}(t) \cdot 1 = 0. \nonumber
\end{eqnarray}
This means that:
\begin{equation} \label{eqn:sum_permanence}
\displaystyle\sum_{s=1}^{M}x_{v,s}(t) = \sum_{s=1}^{M}x_{v,s}(0)~~ \forall t > 0, ~~\forall v \in \mathcal{V}. 
\end{equation}

Imposing that $\boldsymbol{x}_{v}(0) \in \Delta_{M}$, the last equation asserts that 
$\sum_{s=1}^{M}x_{v,s}(t) = 1$ for all time $t > 0$. Joining the results
provided by \eqref{eqn:nonnegativity} and \eqref{eqn:sum_permanence}, we
conclude the following:
\begin{equation} \label{eqn:delta_permanence}
\displaystyle \forall v \in \mathcal{V} : \boldsymbol{x}_{v}(0) \in \Delta_{M} \Rightarrow \boldsymbol{x}_{v}(t) \in \Delta_{M} ~~ \forall t > 0.
\end{equation}
In other words, all trajectories that start inside $\Delta_{M}$ remain inside
$\Delta_{M}$ itself for all time $t > 0$. 
At any time, $\boldsymbol{x}_{v}(t)$ can be
always interpreted as a distribution of strategies.

\subsection{Nash equilibria are rest points of the replicator equation on graph}

Recall that the best response function for the static $(N-M)$-game is:
\begin{equation}\nonumber
\begin{array}{l}
\displaystyle \beta_{v}(\boldsymbol{x}_{-v}) =\\
\displaystyle = \left\{\boldsymbol{x}_{v} \in \Delta_{M} : \pi_{v}(\boldsymbol{x}_{v}, \boldsymbol{x}_{-v}) \geq \pi_{v}(\boldsymbol{z}, \boldsymbol{x}_{-v}) ~~ \forall \boldsymbol{z} \in \Delta_{M}  \right\}.
\end{array}
\end{equation}
Suppose that $\boldsymbol{x}^{*}_{1}(t), \ldots \boldsymbol{x}^{*}_{N}(t)$ is a Nash equilibrium. Then:
\begin{equation} \nonumber
\displaystyle \boldsymbol{x}^{*}_{v}(t) \in \beta_{v}(\boldsymbol{x}^{*}_{-v}(t))
\end{equation}
for each vertex $v$. This means that:
\begin{eqnarray}
\displaystyle \dot{x}^{*}_{v,s}(t) & = & \displaystyle x^{*}_{v,s}(t)(p^{*}_{v,s}(t) - \phi^{*}_{v}(t)) = \nonumber\\
& = & \displaystyle x^{*}_{v,s}(t)(\pi_v(\boldsymbol{e}_{s}, \boldsymbol{x}^{*}_{-v}(t)) - \pi_{v}(\boldsymbol{x}^{*}_{v}(t), \boldsymbol{x}^{*}_{-v}(t))) \leq 0 \nonumber
\end{eqnarray}
Moreover, from \eqref{eqn:delta_permanence} we know that:
\begin{equation}\nonumber
\displaystyle \sum_{s=1}^{M}\dot{x}^{*}_{v,s}(t) = 0 ~~ \forall v \in \mathcal{V},
\end{equation}
and then: 
\begin{equation}\nonumber
\displaystyle \dot{x}^{*}_{v,s} = 0 ~~ \forall v \in \mathcal{V}, ~~\forall s \in \mathcal{S}.
\end{equation}
We can conclude that every Nash equilibrium is also a rest point
of the replicator equation on graph.

\subsection{Pure strategies are rest points of the replicator equation on graph}

Suppose that $\boldsymbol{x}_{v}(t) = \boldsymbol{e}_{q}$. Then $p_{v, q}(t) = \phi_{v}(t)$ and
\begin{equation}\nonumber
\dot{x}_{v,q}(t) = x_{v,q}(t)(p_{v,q}(t) - \phi_{v}(t)) = 1 \cdot 0 = 0.
\end{equation}
In addition, $x_{v,r}(t) = 0$ if $r \neq q$, and again:
\begin{equation}\nonumber
\dot{x}_{v,r}(t) = x_{v,r}(t)(p_{v,r}(t) - \phi_{v}(t)) = 0 \cdot( p_{v,r}(t) - \phi_{v}(t) ) = 0.
\end{equation}
 For this reason:
\begin{equation}\nonumber
\displaystyle \boldsymbol{x}_{v}(t) = \boldsymbol{e}_{q} \Rightarrow \dot{x}_{v,s}(t) = 0 ~~ \forall s,q \in \mathcal{S}.
\end{equation}
This implies that if each $\boldsymbol{x}_{v}(t)$ represents a pure strategy (i.e. 
it is equal to a versor of $\mathbb{R}^{M}$), then we have a 
rest point of the replicator equation on graph.

\subsection{The classical replicator equation as a special case}\label{sec:classical}

Suppose to fix a time lag $\tau > 0$, and assume that the mixed 
strategies are all the same for each vertex and for any time $t_{0} \in [0, \tau)$.
That is:
\begin{equation}\nonumber
\displaystyle \boldsymbol{x}_{v}(t_{0}) = \boldsymbol{c} = [c_{1} \ldots c_{M}]^{T} \in \Delta_{M} ~~ \forall v \in \mathcal{V}, ~~\forall t_{0} \in [0, \tau).
\end{equation}

Consider the payoff model WA and suppose that $\boldsymbol{B}_{v} = \boldsymbol{B}$ for all vertices
$v$.
Following equations \eqref{eqn:pvs} and \eqref{eqn:phiv}, we obtain that
$p_{v,s}(t_{0}) = \boldsymbol{e}_{s}^{T}\boldsymbol{B}\boldsymbol{c}$ and $\phi_{v}(t_{0}) = \boldsymbol{c}^{T}\boldsymbol{B}\boldsymbol{c}$.
In this case, we can rewrite the difference equation 
\eqref{eqn:x_map} as follows:
\begin{equation}
\displaystyle x_{v,s}(t_{0} + \tau) = \frac{c_s(1  + \tau \boldsymbol{e}_{s}^{T}\boldsymbol{B}\boldsymbol{c})}{1 + \tau \boldsymbol{c}^{T}\boldsymbol{B}\boldsymbol{c}}. \nonumber
\end{equation}
Since previous equations do not depend on $v$, we are able to impose that 
$\boldsymbol{y}(t_{0}) = \boldsymbol{x}_{v}(t_{0})$ and $\boldsymbol{y}(t_{0} + \tau) = \boldsymbol{x}_{v}(t_{0} + \tau)$, $\forall v \in \mathcal{V}, ~~\forall t_{0} \in [0, \tau)$,{}
and hence:
\begin{equation}
\displaystyle y_{s}(t_{0} + \tau) = \frac{y_s(t_{0})(1  + \tau \boldsymbol{e}_{s}^{T}\boldsymbol{B}\boldsymbol{y}(t_{0}))}{1 + \tau \boldsymbol{y}(t_{0})^{T}\boldsymbol{B}\boldsymbol{y}(t_{0})}. \nonumber
\end{equation}
It's straightforward to note that any other iteration of the previous map leads to
quantities that are independent from $v$. For example, applying a second iteration
we get that:
\begin{equation}
\displaystyle x_{v,s}(t_{0} + 2\tau) = \frac{y_s(t_{0} + \tau)(1  + \tau \boldsymbol{e}_{s}^{T}\boldsymbol{B}\boldsymbol{y}(t_{0}  + \tau))}{1 + \tau \boldsymbol{y}(t_{0} + \tau)^{T}\boldsymbol{B}\boldsymbol{y}(t_{0} + \tau)}, \nonumber
\end{equation}
and hence we can pose that $\boldsymbol{x}_{v}(t_{0} + 2\tau) = \boldsymbol{y}(t_{0} + 2\tau)$. 
Generalizing to any time lag, $\boldsymbol{x}_{v}(t_{0} + k\tau) = \boldsymbol{y}(t_{0} + k\tau)$ for any non negative
integer $k$. Similarly, $p_{v,s}(t_{0} + k\tau) = \boldsymbol{e}_{s}^{T}\boldsymbol{B}\boldsymbol{y}(t_{0} + k\tau)$ and 
$\phi_{v}(t_{0} + k\tau) = \boldsymbol{y}(t_{0} + k\tau)^{T}\boldsymbol{B}\boldsymbol{y}(t_{0} + k\tau)$ are also
independent from $v$. For these reasons, we pose that 
$p_{s}(t_{0} + k\tau) =p_{v,s}(t_{0} + k\tau)$ and 
$\phi(t_{0} + k\tau) = \phi_{v}(t_{0} + k\tau)$. Then, the
discrete map becomes the following:
\begin{equation}
\displaystyle y_{s}(t_{0} + (k+1)\tau) = \frac{y_s(t_{0} + k\tau)(1  + \tau p_{s}(t_{0} + k\tau))}{1 + \tau \phi(t_{0} + k\tau)}. \label{eqn:y_map_raw}
\end{equation}

Note that, for any $t \geq 0$ there exist a non-negative integer $k$ and
a real number $t_{0} \in [0, \tau)$, with $\tau > 0$ fixed, such that $t = t_{0} + k\tau$.
Then, equation \eqref{eqn:y_map_raw} becomes:
\begin{equation}
\displaystyle y_{s}(t + \tau) = \frac{y_s(t)(1  + \tau p_{s}(t))}{1 + \tau \phi(t)} ~~ \forall t \geq 0. \label{eqn:y_map} 
\end{equation}
Considering the difference ratio $\frac{1}{\tau}(y_{s}(t + \tau) - y_{s}(t))$,
and letting $\tau \rightarrow 0$, we obtain the following differential equation:
\begin{equation}\label{eqn:re_classic}
\displaystyle \dot{y}_{s}(t)  = y_{s}(t)(p_{s}(t) - \phi(t)),
\end{equation}
which is the classical replicator equation.\\

This result is quite straightforward if we imagine to divide a wide population 
of replicators into $N$ subpopulations, assuming that all of them are described by
the same mixed strategy of the total one at initial time. Then, each 
subpopulation will behave exactly as the total one. Hence, the dynamics
of a single subpopulation in a vertex can be described 
by the classical replicator equation applied to the single population,
whatever is the graph used.

 
\section{Simulations}\label{sec:sim}

\begin{figure}[h!]
\begin{center}
\includegraphics[width=.6\columnwidth]{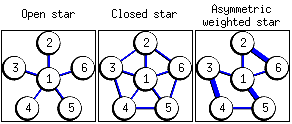}
\caption{\footnotesize{Graphs topologies. All edges have weights equal to $1$, except for
thicker ones in the asymmetric weighted graph, having weights equal to $3$.
}}
\label{fig:graphs}
\end{center}
\end{figure}

In this chapter, we present some simulations produced by
equation \eqref{eqn:ode}. The $WA$ payoff
model is used. In particular,
we set up experimental sessions by considering different $2$-strategies payoff 
matrices ($\mathcal{S} = \left\{1, 2\right\}$); it is assumed that every 
vertex has the
same payoff matrix. Each session has been developed over
$3$ different graphs with $6$ vertices as reported in Figure \ref{fig:graphs}.

All edges represented in Figure \ref{fig:graphs} have the same weight, except for
thicker ones in the asymmetric weighted graph. Note that we are using only undirected graphs
(i.e. $a_{v,w} = a_{w,v}$).
Our aim is to show the behavior of the replicator equation on graphs
when initial players strategies are almost pure. In fact, a vertex player with a 
pure strategy is in steady state; for this reason, initial 
conditions used for vertex players are equals to slightly perturbed 
pure strategies (i.e. $[0.99 ~ 0.01]^{T}$ and $[0.01 ~ 0.99]^{T}$ 
are used in place of pure strategies $1$ and $2$, respectively). 
Replicator equation on graphs has been simulated until a steady
state behavior is reached, starting from $4$ different distribution initial conditions on the graph.\\

The steady state situations are shown in Figures \ref{fig:sym_bistable}, \ref{fig:asym_bistable},
\ref{fig:prisoner} and \ref{fig:coexistence}. The first
column of each Figure gives a picture of the initial conditions used,
while others report the solution of the simulations when steady state
is reached for each of the considered graphs. The color of each vertex
indicates the value
of $x_{v,1}$, and hence it visually quantifies the inclination of player
$v$ toward one of the $2$ feasible pure strategies;
yellow is used for player with strategy $1$ ($x_{v,1} = 1$), red is 
for strategy $2$ ($x_{v,1} = 0$). Mixed strategies ($0<x_{v,1}<1$) are 
indicated by shaded colors, according to the color bar
at the bottom of the Figures. Moreover, Figures \ref{fig:timecourse_pd}
and \ref{fig:timecourse_coex} report the dynamical evolution obtained
on the asymmetric weighted graph; the same initial condition is used 
in both Figures, while payoff matrices are different.
The following sections will discuss in detail the results of each simulation.

\subsection{Two pure Nash equilibria}

\begin{figure}[h!]
\begin{center}
\includegraphics[width=0.8\columnwidth]{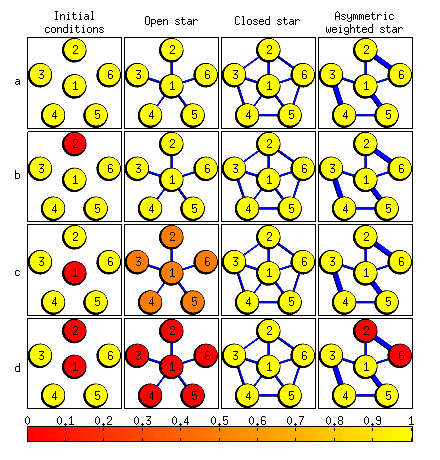}
\caption{\footnotesize{Strategies distribution at $t=0$ (first column) and at $t=50$ (other columns) on open, closed and 
asymmetric weighted star graphs, when the payoff matrix is reported in equation \eqref{eqn:twopure}
and $\theta = 1$. $4$ different initial conditions are 
considered: homogeneous (a), external outlayer (b), central outlayer 
(c) and external-central outlayers (d).
}}
\label{fig:sym_bistable}
\end{center}
\end{figure}

In this first experimental session, we used the following payoff matrix:
\begin{equation}\label{eqn:twopure}
\boldsymbol{B} = \left[\begin{array}{cc}
1 & 0\\0 & \theta
\end{array}\right],
\end{equation}
with $\theta > 0$.

The $2$-players game described by $\boldsymbol{B}$ has $2$ strict pure Nash equilibria
(i.e. both players use strategy $1$ or $2$) and a mixed Nash equilibrium
$$\displaystyle\boldsymbol{x}^{*} = \left[\frac{\theta}{1+\theta} ~~ \frac{1}{1 + \theta}\right]^{T}.$$
The classical replicator equation, 
based on matrix $\boldsymbol{B}$, has exactly $3$ rest points which coincide
with the Nash equilibria reported above. Moreover, mixed equilibrium is 
repulsive, while pure equilibria are attractive; for this reason, we say
that $B$ is a bistable payoff matrix. \\

Figure \ref{fig:sym_bistable} reports some results obtained when $\theta = 1$.
Row (a) of the Figure shows what happens
when an homogeneous initial condition is used; as said in section \ref{sec:classical},
the dynamics is the same for each vertex player, and it is equivalent
to the solution given by the classical replicator equation, whichever is
the underlying graph structure.
After a certain time, all vertex players adopt pure strategy $1$, since
it represents an attractive rest point, and initial condition is
in inside the relative basin of attraction.\\

In the row (b) of Figure \ref{fig:sym_bistable} are reported the
steady state situations obtained by using an homogeneous initial
condition,
where only one peripheral player uses the quasi-pure strategy $2$. At the
end of simulation, the pure strategy $1$ spreads all over the considered graphs.
Let's consider the open graph situation: the vertex player $1$, 
which is the unique neighbors of player $2$, has no
will to change his own strategy, since he is surrounded by $4$ yellow players.
Similarly, on the closed and asymmetric weighted star, neighbors
of player $2$ see an equivalent player which is almost yellow. Thus, none
of them wants to change, and player $2$ must modify his strategy.
Hence player $2$ must change his strategy to obtain a good payoff.
In a certain way, the ''rebel`` peripheral player decides to adapt himself
to the majority.\\

The dynamical behavior is slightly different when the central hub is the rebel. 
In the row (c) of the Figure \ref{fig:sym_bistable} are shown the solutions
of the replicator equation on graphs for this initial condition.
When the open star is used, player $1$ sees a yellow equivalent player, while all peripheral players
have only him as neighbor. Player $1$ decides to change his own strategy to yellow,
while all others do the exact opposite. After a certain time,
they meet half way, at the mixed equilibrium $[0.5 ~ 0.5]^{T}$. The different
position of the rebel player in the graph influences a lot the dynamics 
of the whole system; the leader (player $1$) understands that he must
modify his own strategy according to his neighborhood, while all other players do
the same, since their only opponent is player $1$ himself. However, closed and asymmetric weighted
graphs are more resistant to the influence of player $1$, because
the peripheral players have more than one neighbors; in these situations,
player $1$ does not play anymore as a leader able to change the whole dynamics.\\

The last row (d) of Figure \ref{fig:sym_bistable} reports the final
solutions when both player $1$ and $2$ use the quasi-pure strategy $2$.
While the closed star structure remains resistant to the influence of 
rebel players, the other graphs do not. The open star becomes all red
at final time. This is because player $2$ sees only player $1$: they are
both red, so player $2$ doesn't want to change strategy. 
Simultaneously, yellow neighbors of $1$ change their strategy to red,
since they see only a red player. \\

\begin{figure}[h!]
\begin{center}
\includegraphics[width=0.8\columnwidth]{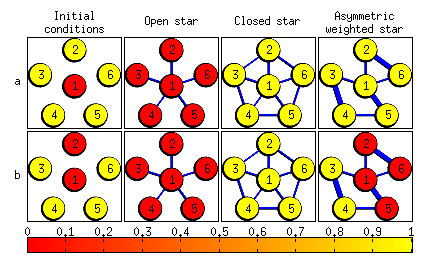}
\caption{\footnotesize{Strategies distribution at $t=0$ (first column) and at $t=50$ (other columns) on open, closed and 
asymmetric weighted star graphs, when the payoff matrix is reported in equation \eqref{eqn:twopure}
and $\theta = 1.1$. $2$ different initial conditions are 
considered: central outlayer (a) and external-central outlayers (b).
}}
\label{fig:asym_bistable}
\end{center}
\end{figure}

Changing the value of the parameter $\theta$ leads to different behaviors.
When $\theta < 1$, first strategy becomes stronger and it spreads all
over the considered graphs as $\theta$ goes to $0$.
In Figure \ref{fig:asym_bistable} are reported some results obtained
with $\theta = 1.1$. In particular, when player $1$ uses strategy $2$
at the beginning, then mixed equilibrium is not reached anymore on the open
star graph; all vertices adopt strategy $2$, which is slightly better than 
strategy $1$. The strength of strategy $2$ is also visible on the asymmetric
weighted star, when at the beginning both player $1$ and $2$ adopt strategy $2$;
in Figure \ref{fig:sym_bistable} ($\theta=1$) we have shown that on
steady state, players $2$ and $6$ are the only ones red, while when $\theta=1.1$,
also players $1$ and $5$ do. In general, when $\theta > 1$, strategy $2$ 
becomes stronger and it spreads all over the considered graphs as $\theta$ grows up.

\subsection{Prisoners' dilemma}

In this section, we show the results obtained with the replicator equation on
graphs, by using a modified version of the prisoners' dilemma game, as proposed 
in \cite{N2006}. The payoff matrix is the following: 
\begin{equation}\label{eqn:prisoner}
\boldsymbol{B} = \left[\begin{array}{cc}
1 & 0\\\theta & 0
\end{array}
\right],
\end{equation}
where $\theta > 1$.\\

\begin{figure}[h!]
\begin{center}
\includegraphics[width=0.8\columnwidth]{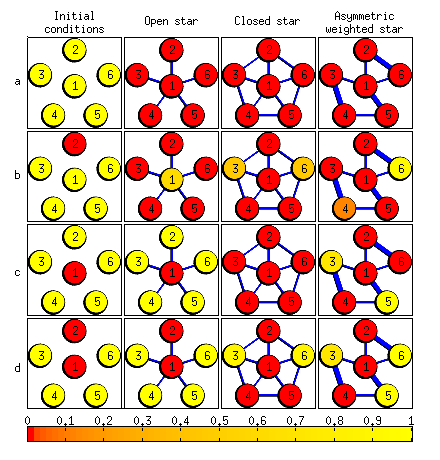}
\caption{\footnotesize{Strategies distribution at $t=0$ (first column) and at $t=100$ (other columns) on open, closed and 
asymmetric weighted star graphs, when the payoff matrix is reported in equation \eqref{eqn:prisoner}
and $\theta = 1.5$. $4$ different initial conditions are 
considered: homogeneous (a), external outlayer (b), central outlayer 
(c) and external-central outlayers (d).
}}
\label{fig:prisoner}
\end{center}
\end{figure}

\begin{figure}[h!]
\begin{center}
\includegraphics[width=0.8\columnwidth]{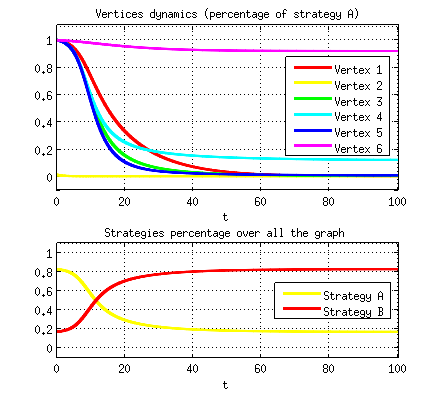}
\caption{\footnotesize{Prisoners' dilemma game. Time courses of each vertex (top) and average behavior
over the network (bottom). Network topology, initial conditions and parameter $\theta$ are the same used
in Figure \ref{fig:prisoner}, row b, column 4.
}}
\label{fig:timecourse_pd}
\end{center}
\end{figure}

\textit{Cooperate} and \textit{Defect} are the names typically used to indicate, respectively,
the strategy $1$ and $2$ of this classic game. 
The dilemma is that mutual cooperation produces a better outcome than 
mutual defection; however, at the individual level, the choice to cooperate
is not rational from a selfish point of view. In other words, the
$2$-players 
game has only one Nash equilibrium, reached when both players defect.
Note that in this version of the prisoner's dilemma, the Nash equilibrium is
non strict. Moreover, classical replicator equation based upon payoff matrix reported in equation
\eqref{eqn:prisoner}, has $2$ rest points, which correspond to the  pure strategies $(1,1)$ and $(2,2)$. In particular, the first one is repulsive,
while the latter is attractive.\\

Although mutual defection represents both a Nash and a dynamical equilibrium,
many works have shown that cooperation does not vanishes when games are
played over graphs and $\theta$ is equal to suitable values
(see \cite{TK2002,N2006,SRP2006,GR2012}). The resilience of cooperation is shown
in \ref{fig:prisoner}, where $\theta$ is set to $1.5$. Steady states depend on the initial conditions and on the type
of graph used,
and behaviors can be very heterogeneous. When an homogeneous initial condition is considered (row (a)),
all players on graphs become defectors (again, this is the case when
the classical and proposed replicator equations are the same). On the other side,
when initial conditions are not homogeneous (rows (b), (c) and (d)), cooperation
does not always completely vanish. Figure \ref{fig:timecourse_pd} shows
the time course of the variable $x_{v,1}$ for each vertex of the graph.
In particular, the initial conditions with external outlayer and the
asymmetric weighted graph have been used. 

\subsection{Unique mixed Nash equilibrium}

In some $2$-players games there are no pure Nash equilibria.
Nevertheless, Nash theorem guarantees that at least a mixed equilibrium exists.
For example, this happens when payoff matrix is defined as follows:
\begin{equation}\label{eqn:coexistence}
\boldsymbol{B} = \left[\begin{array}{cc}
0 & 1\\1 & 0
\end{array}
\right].
\end{equation}
The unique mixed Nash equilibrium is:
$$\displaystyle\boldsymbol{x}^{*} = \left[0.5 ~~ 0.5\right]^{T}.$$
Classical replicator equation has $3$ rest points; symmetric couples
of pure strategies, $(1,1)$ and $(2,2)$ are repulsive, while the 
mixed equilibrium is attractive. In this case, we speak about $coexistence$ of both feasible strategies.\\

\begin{figure}[h!]
\begin{center}
\includegraphics[width=0.8\columnwidth]{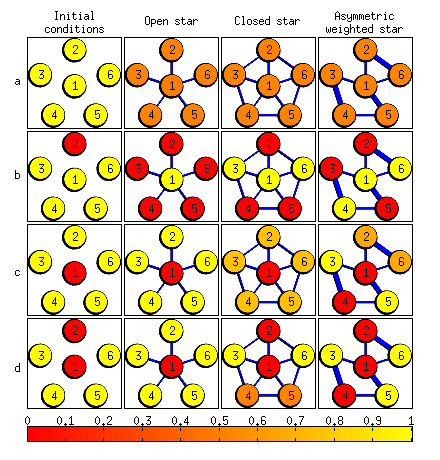}
\caption{\footnotesize{Strategies distribution at $t=0$ (first column) and at $t=50$ (other columns) on open, closed and 
asymmetric weighted star graphs, when the payoff matrix is reported in equation \eqref{eqn:coexistence}. $4$ different initial conditions are 
considered: homogeneous (a), external outlayer (b), central outlayer 
(c) and external-central outlayers (d).
}}
\label{fig:coexistence}
\end{center}
\end{figure}

\begin{figure}[h!]
\begin{center}
\includegraphics[width=0.8\columnwidth]{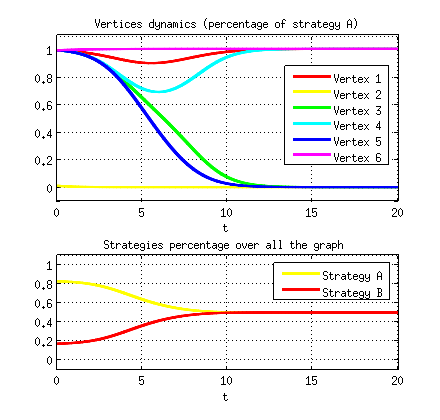}
\caption{\footnotesize{Unique mixed Nash equilibrium game. Time courses of each vertex (top) and average behavior
over the network (bottom). Network topology and initial conditions are the same used
in Figure \ref{fig:coexistence}, row b, column 4.
}}
\label{fig:timecourse_coex}
\end{center}
\end{figure}

In Figure \ref{fig:coexistence} the steady state
solutions when payoff matrix defined in equation \eqref{eqn:coexistence} is used, are reported.
Again, when we have an homogeneous initial condition, everything
works like a classical replicator equation, and hence, all players
go to the mixed Nash equilibrium. When initial condition is not
homogeneous, behaviors obtained through the replicator equation on 
graphs are strongly based on the topological structure of the 
underlying graph and on the initial conditions.
Figure \ref{fig:timecourse_coex} shows in details the behavior of
the population when the asymmetric weighted graph  and initial conditions
with external outlayer are supposed.

 
\section{Conclusions}\label{sec:conclusions}

In this work a new mathematical model for evolutionary games on graphs
with generic topology has been developed. We proposed a replicator 
equation on graphs,dealing with a finite population of players
connected through an arbitrary topology.
A link between two players can be weighted by a positive
real number to indicate the strength of the connection. Furthermore, 
the different perception that each player has about the game is 
modeled by allowing the presence of directed links and different
payoff matrices for each member of the population. A player obtains his
outcome after $2$-players games are played with his neighbors; payoffs
of each game are averaged (WA model) or simply summed up (WS model). 
Moreover, it has been shown that the proposed replicator equation on 
graphs extends the classical one, under the hypotheses that WA model for payoffs
is used, homogeneous initial conditions 
over the vertices are considered, all vertex players have the same
payoff matrix. In any case, no limitations are imposed to the underlying graph.\\

Experimental results showed that the dynamics of evolutionary games 
are strongly influenced by the network topology. As expected, more complex
behavior emerges with respect to the classical replicator equation. For example,
in the prisoner's dilemma game, cooperative and non-cooperative behaviors
can coexist over the graph. Moreover, when a $2$-player game with strictly 
dominant strategies is considered, heterogeneous behavior is obtained, i.e.
a part of the population chooses to play a dominant strategy, while others
use different strategies. Then, players become mixed (coexistence of strategies).\\

The very first step for extending this work is the study of dynamical
and evolutionary stability of the rest points. By the way, we imagine that
the concept of evolutionary stability must be revisited to deal with the
proposed evolutionary multi-players game model based on graph, for which a theoretical effort
is needed. Indeed, in our opinion, the basic question
``is strategy $s$ resistant to invasion?'' must be reformulated to fit with the new model, 
where the population of players is finite and is organized according to a social structure.\\

The theory developed in this paper can also be extended to 3 or more strategies and can consider 
more complex topologies of the graph, such as small world, scale free, and random complex networks.\\
From an applicative point of view, the authors intend to use the replicator equation
on graphs to deal with biological and physical processes, 
such as bacterial growth \cite{BBGMMPV2012}, model of
brain dynamics \cite{MCSM2013} and reaction-diffusion phenomena \cite{MMS2011}.  
The developed model can be also profitably applied to solve networked socio-economics problems, 
such as decision making for the development of marketing strategies.



\end{document}